# Statistical properties of the phase transitions in a spin model for market microstructure

Muffasir Badshah
Department of Mathematics, Drexel University
3141 Chestnut Street
Philadelphia, PA 19104
+1-215-895-2680

mhb25@drexel.edu

Dr. Robert Boyer
Department of Mathematics, Drexel University
Korman Center, Room 264
Philadelphia, PA 19104
+1-215-895-1854

rboyer@drexel.edu

Dr. Ted Theodosopoulos
Department of Decision Sciences and Department of Mathematics, Drexel University
Academic Building, Room 224
Philadelphia, PA 19104
+1-215-895-6997

theo@drexel.edu

## ABSTRACT
Increased day-trading activity and the subsequent jump in intraday volatility and trading volume fluctuations has raised consideratble interest in models for financial market microstructure. We investigate the random transitions between two phases of an agent-based spin market model on a random network. The objective of the agents is to balance their desire to belong to the global minority and simultaneously to the local majority. We show that transitions between the "ordered" and "disordered" phases follow a Poisson process with a rate that is a monotonically decreasing function of the network connectivity.

## Categories and Subject Descriptors
G.3 [**Mathematics of Computing**]: Probability and Statistics – *probabilistic algorithms and statistical computing*.

J.3 [**Social and Behavioral Sciences**]: Economics.

## 1. MODEL OVERVIEW
In this paper we investigate a spin market model introduced Bornholdt and collaborators in [1,2]. The model consists of agents on a network, with vertex set *V* and edge set *E*. Each agent's state is a binary variable $S_i(t)$ (spin), where *i* refers to the agent. At time *t+1*, agent *j* sets their spin to +1 (or -1) with probability $p_+(t) = (1 + \exp\{-2\beta h_i(t)\})$ (or $p_-(t) = 1 - p_+(t)$), where the Hamiltonian, $h_i$ is given by

$$h_i(t) = \sum_{(i,j) \in E} S_j(t) - \alpha S_i(t) \frac{1}{N} \left| \sum_{j=1}^{N} S_j(t) \right|$$

Agents in this model try to balance the desire to be in the global minority with the desire to be in the local majority. The former effect is motivated by the minority game, which captures the contrarian tendencies in the market. Clearly if there is a global imbalance in favor of buys, you'd rather be a seller. On the other hand the local majority effect is reminiscent of the voter model, in which agents learn from the opinions of their neighbors.

This model was studied under an asynchronous updating rule in [3] and its asymptotic equilibria were characterized. Here we investigate its dynamics using synchronous updates. In this setting, it is known that, for appropriately chosen values of $\alpha$ and $\beta$, the system can be found in one of two phases: the traders can be switching continuously between the two spins ("disordered" phase), or they may be segregated into stable coalitions of buyers and sellers ("ordered" phase).

The network we used in this paper consists of sixteen agents, one at each vertex. Three basic neighborhood structures were investigated, together with random variants. Specifically we analyzed the 2-neighbor, 4-neighbor and 8-neighbor lattices shown in Figures 1-3. Figure 4 shows the labeling of the neighbors that is used in describing the perturbed networks that were also studied.

In the simulations that follow we use the values $\alpha = 4$ and $\beta = 0.5$. Bornholdt has shown that this stochastic model exhibits two phases, an "ordered" one in which the Hamiltonian of a particular site stays within a narrow band around 0, and a "disordered" on, characterized by large fluctuations of the Hamiltonian. In order to extract and analyze the intervals of the ordered regime (stable behavior of $h_i(t)$), the following conditions were mandated: $|h_i(t+1) - h_i(t)| < 0.5$. A sample of this extraction is shown in Figure 5.

In each case, simulations of 8192 steps were performed and the time periods spent in the "ordered" vs. "disordered" phase were recorded.

## 2. SIMULATION RESULTS

The first observation we made was that the overall time spent in the "ordered" phase is an increasing function of the connectivity of the network. This can be seen in Figure 6b, which shows that the 8-neighbor model spends over 95% of the time in the "ordered" phase, as compared to less than 12% for the two-neighbor model. The mean ratios for the three models are provided in Figure 6b. Single factor ANOVA (p-Value = 0) with Tukey-Kramer procedure were used to discern differences among the mean ratios of the three groups with 95% certainty.

When we begin eliminating links from the 8-neighbor model, we see in Figure 7b that the resulting percentages follow the same pattern. For instance eliminating two random links from the neighborhood of each agent leaves a randomly oriented set of six edges incident to each vertex, and leads to 70% of the simulation time spent in the "ordered" phase. This percentage increases to 86% when only one random link is eliminated.

Also, we investigated the duration $T$ of each period spent in the "ordered" phase. The graphs of $log(Pr(T>t))$ versus time t are given for each described category and tables are provided for the ratio of overall time in the ordered regime to the total time. The distribution suggests the following exponential decay behavior: for the two-neighbor model $\Pr(T > t) \sim e^{-0.44t}$, for the four-neighbor model $\Pr(T > t) \sim e^{-0.073t}$ and for the eight-neighbor model $\Pr(T > t) \sim e^{-0.016t}$. Single factor ANOVA (p-Value = 2.57E-13) with Tukey-Kramer procedure were used to discern differences among the means of the three groups with 95% certainty. The means and plots are provided in the Figure 6a for all three models.

The results in the model 2, 4 and 8 suggest strongly that as the network connectivity increases, the model becomes more stable. To further explore the geometric structure dependencies, we analyze the 8-neighbor model more extensively. Specifically, we proceed to eliminate random links from the eight-neighbor lattice, thus creating random networks with constant degree distributions.

Figure 7a shows the resulting distribution of the duration of each period spent in the "ordered" phase: eliminating 1 random link (7 links present) leads to $\Pr(T > t) \sim e^{-0.023t}$, eliminating 2 random links (6 links present) leads to $\Pr(T > t) \sim e^{-0.144t}$ while finally eliminating 3 random links (5 links present) leads to $\Pr(T > t) \sim e^{-0.3797t}$. Single factor ANOVA was used to discern differences among groups (p-Value = 3.16E-17). Tukey-Kramer procedure determined that there exists a difference among the mean slopes of all three groups with 95% certainty.

On the other hand, Figure 7b shows the overall percentage of simulation steps spent in the "ordered" phase, when we eliminate random links. Single factor ANOVA (p-Value = 0) with Tukey-Kramer procedure determined that there exists a difference among the mean ratios of all three groups with 95% certainty.

We observe that the percentage of time spent in the "ordered" phase is a strongly increasing function of the network connectivity. This effect is captured by Figure 8. In particular, we see that this relationship is power law, with exponent around 1.5.

Figure 9 on the other hand shows the monotonically decreasing relationship between the Poisson rate for the transitions between the "ordered" and "disordered" phases and the network degree. This relationship is also a power law. However, in this case we obtain two different exponents depending on whether the network is a lattice (lower branch of the graph) or a random network (i.e. with random missing links).

## 3. FIGURES/CAPTIONS

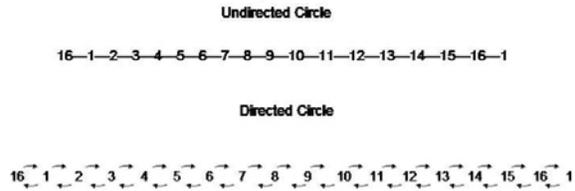

**Figure 1. Two Neighbors – Model 2**

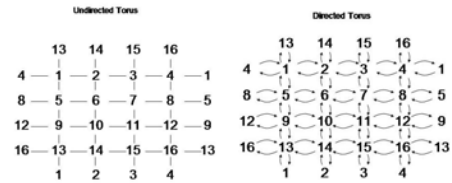

**Figure 2. Four Neighbors – Model 4**

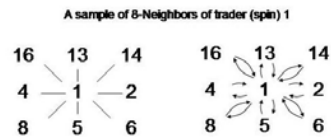

**Figure 3. Eight Neighbors – Model 8**

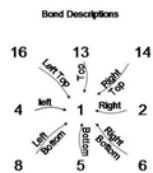

**Figure 4. Bond Classification**

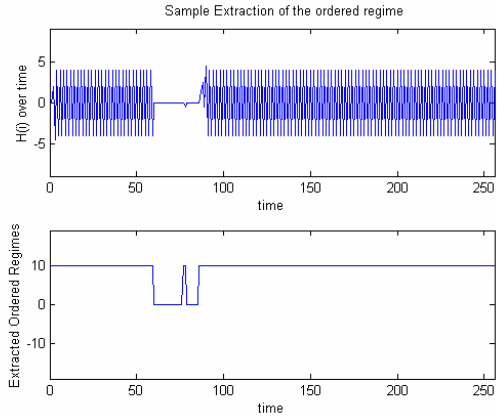

**Figure 5. Sample Extraction of the Ordered Regime**

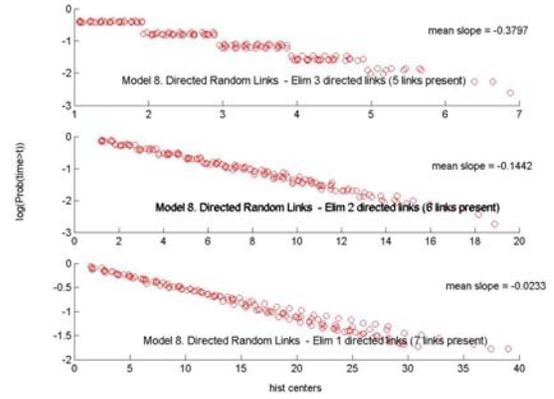

**Figure 7a. Model 8 – Randomly eliminated one, two and three directed links. Analysis of duration of time Interval/s in the Ordered Regime**

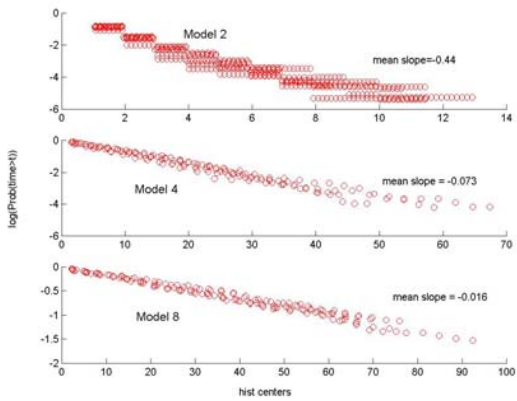

**Figure 6a. Model 2, 4 and 8 - Analysis of the duration of time Interval/s in the Ordered Regime**

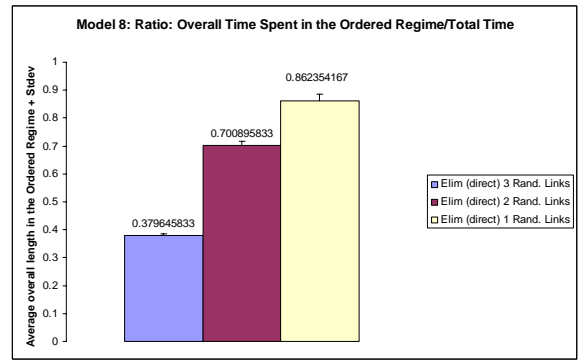

**Figure 7b. Model 8 - Randomly eliminated one, two and three directed links. Analysis of the ratio – total duration in the ordered regime versus overall length of time**

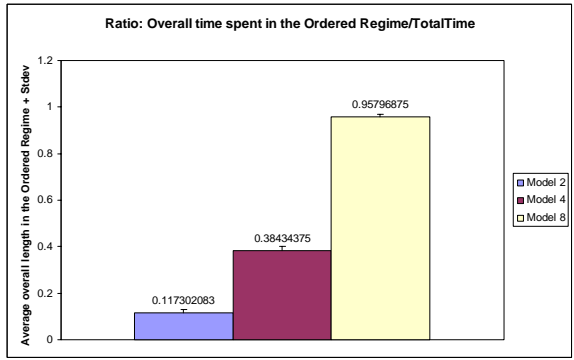

**Figure 6b. Model 2, 4 and 8 - Analysis of the ratio – total duration in the ordered regime versus overall length of time**

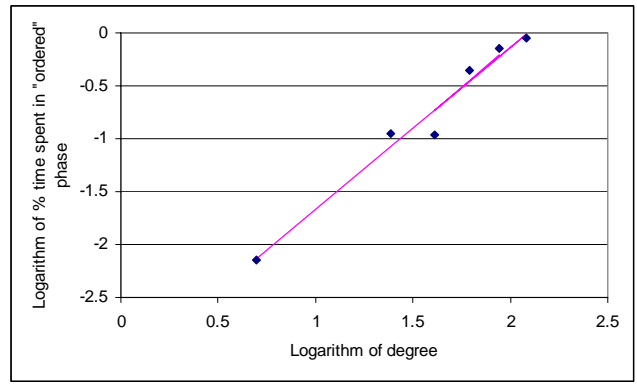

**Figure 8. Power law dependence of percentage of time spent in the "ordered" phase as a function of the degree of the network.**

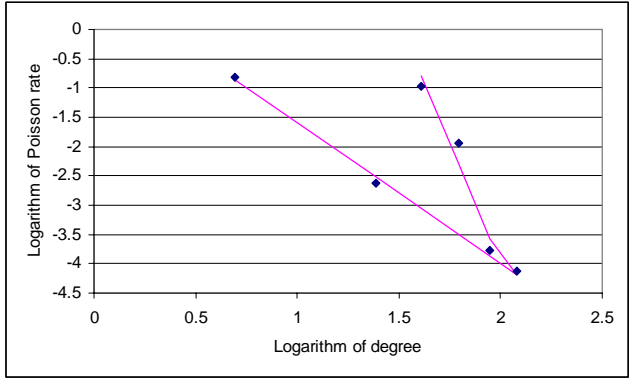

**Figure 9. Power law dependence of Poisson rate as a function of the degree of the network. The steeper curve corresponds to the randomly perturb networks.**

## 4. CONCLUSION AND FUTURE WORK

We analyzed the transitions of an agent-based spin market model between its two phases. Our analysis indicates that the time spent in the "ordered" phase as a function of the network connectivity is an increasing power law, which appears to be robust to random perturbations of the network structure.

Also, we conclude that the transitions from the "ordered" phase to the "disordered" phase follow a Poisson process. The rate of this process exhibits a decreasing power law dependence on the connectivity of the underlying network, which appears to be quite sensitive to random perturbations of the fine network structure.

The reduction in frequency and stability of the ordered regime can be explained by the new traders joining in the majority group [1,2]. The frequency of risk averters leaving the majority coalition is less than the incoming traders joining the coalition. An alternative explanation to the stability of trader behavior with the increase in size of the local neighborhood can be derived from herd behavior [4,5].

In a follow-on paper we will examine the impact of the phase transitions studied here to market volatilityand volume fluctuations. Also, a natural extension of the work presented here is the investigation of random networks with non-constant degree distributions.

## 5. ACKNOWLEDGMENTS

Our thanks to ACM SIGCHI for allowing us to modify templates they had developed.

We would also like to thank the Lebow College of Business and the Department of Mathematics at Drexel University for their facilities, equipment and most importantly, their support.

## 6. REFERENCES

[1] Bornholdt, S., *Expectation Bubbles in a spin model of markets: intermittency from frustration across scales.* Int. J. of Modern Physics, Vol. 12, 5 (2001), 667-674.

[2] Bornholdt, S., Fujiwara, Y., and Taisei, K. *Dynamics of price and trading volume in a spin model of stock markets with heterogeneous agents.* Physica A, 316 (2002), 441-452.

[3] Theodosopoulos, T. *Uncertainty relations in models of market microstructure.* www.arXiv.org, Sep. 2004.

[4] Cont, R., and Bouchaud, J. *Herd Behavior and aggregate fluctuations in financial markets.* Macroeconomic Dynamics, 4 (2000), 170-196.

[5] Ghoulmie, F., *Switching by agents between two trading behaviors and the stylized facts of financial markets.* www.arXiv.org, Aug. 2004.